\documentclass[preprint,12pt]{elsarticle}

\usepackage{amssymb}

\usepackage{epsfig}
\usepackage{amsmath}

\newtheorem{thm}{Theorem}

\newtheorem{lemma}[thm]{Lemma} 
 
\newtheorem{cor}[thm]{Corollary} 
 
\newtheorem{exa}[thm]{Example}

\newtheorem{defn}[thm]{Definition}

\newcommand{\beq}{\begin{eqnarray}}
\newcommand{\eeq}{\end{eqnarray}}
\newcommand{\V}{\mathbb V}

\journal{Discrete Mathematics}

\begin{document}

\begin{frontmatter}



\title{Rectangularity}


\author{Tim Boykett}

\address{Institute for Algebra, Johannes Kepler Universita\"at Linz \\ and Time's Up Research}
\ead{tim@timesup.org}
\ead[url]{http://algebra.uni-linz.ac.at  http://timesup.org}

\begin{abstract}
We introduce a condition on arrays in some way maximally distinct from Latin square condition, as well as some other conditions on algebras, graphs and $0,1$-matrices. We show that these are essentially the same structures, generalising a similar collection of models presented by Knuth in 1970.

We find ways in which these structures can be made more specific, relating to existing investigations, then show that they are also extremely general; the groupoids satisfy no nontrivial equations. Some construction methods are presented and some conjectures made as to how certain structures are preserved by these constructions.
Finally we investigate to what degree partial arrays satisfying our conditions and partial Latin squares overlap.

Note that this paper is slightly updated from the first submission.

\end{abstract}


\end{frontmatter}

\section{Introduction}

Latin squares have been used and analysed for centuries. 
They are of great interest in themselves as well as for their connections to a number of other areas in combinatorics and algebra. From a Latin square one obtains immediately an algebraic structure known as a quasigroup, much work has investigated these and related objects such as loops.
In this paper we consider a class of structures that are somehow maximally unlike Latin squares, but use some similar ideas to approach them. We will use graph theory, combinatorics and algebra in order to investigate the properties of these structures.

The investigation uses a spectrum of approaches to understanding the structures of interest. While the origins derive from earlier work, the first part of this connecting work to appear was Knuth's \cite{knuth70} equivalences between a graph theory problem, a matrix formulation \cite{hoffman67} and an algebraic structure discussed by Trevor Evans\cite{evans67a}. These have been investigated at length by a number of researchers since then (see details in Section \ref{seccentral}). The author came across related structures with applications in computer science \cite{boykett04}. The current work arises from a further  generalisation of the two areas of investigation.

We will look at four distinct models and show that the structures are intimately related. We will then put the various special classes of structures into a relationship with one another. We will see that the class of structures is extremely general, the algebras lying within no nontrivial variety. Examples can be constructed using partitions of a point set such as with group factorisations and examples can be combined in a number of ways. Finally we look at the common partial structures between our rectangular ones and Latin squares.

\section{Models and Motivations}

In this section we introduce several models from combinatorics and algebra, 
before showing that these are equivalent.

In a Latin square of order $n$ we have every every row and every column containing precisely one copy of each number in $\{1,\ldots,n\}$. 
One can equivalently state that in every row and every column, each pair of elements appears. We have reached the maximum of getting as many pairs of elements in each row and column. The converse question arises: how can we fill an array so that the lowest number of pairs occurs in each row and column?
It turns out we can do this and require that the rows and columns are  pair disjoint, i.e.\ if a pair appears in a row, it never appears in a column.
Here we define an array that is in some sense maximally unlike a Latin square.
\begin{defn}
Let $M$ be an $n\times n$ array with entries from $\{1,\ldots,n\}$. $M$ has property $P_1$ iff when two elements appear in one row together, they never appear in one column together and vice versa.
\end{defn}

A simple example is to fill the array entirely with one element. Then no pairs occur in rows or columns and we are trivially finished. 
We will call an array \emph{full} if all elements in  $\{1,\ldots,n\}$ arise in the array.
The following two arrays satisfy $P_1$. Note that the second is maximally pair disjoint: every pair occurs in some row or some column, which is not the case in the first example. We will call such arrays \emph{maximal}.

\[
\begin{array}{ c c c c }
 1& 1& 3 & 3  \\
 2 & 2 & 4 & 4 \\
 1 & 1 & 3 & 3 \\
2& 2& 4 & 4 
\end{array}
\hspace{12mm}
\begin{array}{ c c c c }
 1 & 1& 4 & 4  \\
 2 & 2 & 3 & 3 \\
 1 & 1 & 3 & 3 \\
 2 & 2 & 4 & 4 
\end{array}
\]

\begin{defn}
Let $N=\{1,\ldots,n\}$, $(N,R)$ and $(N,G)$ be two graphs on the node set $N$ that we will call the red and green graphs. We say this graph pair has property $P_2$ if for every pair of nodes $a,b\in N$ there is a unique red-green path, i.e.\  $\exists ! c \in N$ s.t.\ $(a,c) \in R$ and $(c,b)\in G$.
\end{defn}

One could talk about these as an idealised product distribution graphs. If every node on $N$ represents a producer and a consumer we use $R$ to represent the transport to a distribution center and $G$ to represent the transport from a distribution center to the consumer. 
For instance a farmers' market as a unique distribution center has a selected 
node $a\in N$ with $R=\{(x,a): a\in N\}$, all farmer's take their produce to the market 
at $a$ and $G=\{(a,x): x \in A\}$, the farmers take what they 
need back from the market to their farms.

\begin{defn}
Let $(A,*)$ be a $(2)$-algebra such that
\[
a*b=c*d=x \Rightarrow a*d=c*b=x
\]
for all $a,b,c,d,x \in A$. We call $A$ a \emph{rectangular groupoid}.
\end{defn}
Rectangular groupoids form a quasivariety as they are defined by an implication \cite{bursank}. We will see below that they form a proper quasivariety, i.e.\ the class of rectangular groupoids is not closed under taking homomorphic images.

Note that the implication 
\beq
a*b=c*d \Rightarrow a*d=a*b \label{exrect}
\eeq
 is sufficient to show rectangularity by the symmetry of the equality relation.

\begin{defn}
Let $A,B$ be two $n\times n$ $0,1$-matrices. We say $A,B$ have the property $P_4$ iff $AB=J$, the matrix consisting of all 1s. 
\end{defn}

We proceed now to show that these four concepts are closely related.
This first result echoes the connection between Latin squares and quasigroups.
\begin{thm}
An $n \times n$ array $M$ has $P_1$ iff it is the Cayley table of a rectangular groupoid $(\{1,\ldots,n\},*)$.
\end{thm}
Proof: $(\Rightarrow)$ Suppose $M$ satisfies $P_1$. Let $A=\{1,\ldots,n\}$ and define $a*b=c$ for $a,b\in A$ with $c$ the $(a,b)$ entry in $M$.
 Now suppose $a*b=c*d$ for some $a,b,c,d \in A$, let $x=a*b=c*d$ and $y=a*d$.  Then both $x$ and $y$ are in the $a$ row and the $d$ column. Thus $x=y$ so $a*d=a*b$ and $(A,*)$ is rectangular.
 
 $(\Leftarrow)$ Let $(A,*)$ be a rectangular groupoid and label $A=\{1,\ldots,n\}$. 
 We create the array $M$ with entry $(a,b)$ equal to $a*b$.
 Suppose two elements $x,y$ appear in some column and in some row. Let the row be $a$ and the column be $b$.
 Then there exist some $c,d\in A$ such that $a*c=x$ and $d*b=x$ so by the rectangularity property $a*b=x$. However the same argument applies to $y$ in the same row and column so $a*b=y$ so $x=y$ and we see that our array satisfies $P_1$.
 \hfill$\Box$

The following result is a direct application of what an incidence matrix means.
\begin{thm}
Two graphs $(N,R)$ and $(N,G)$ have property $P_2$ iff their node-node incidence matrices $I_R,I_G$ have property $P_4$. 
\end{thm}
Proof: The $(i,j)$ entry in the product $I_RI_G$ counts how many length 2 paths from node $i$ to node $j$ exist with the first edge in $(N,R)$ and the second edge in $(N,G)$. Thus the graph pair $(N,R),(N,G)$ satisfies $P_2$ iff 
$I_RI_G$ has a 1 in each entry iff $I_R,I_G$ have property $P_4$. 
\hfill$\Box$

The following two results bind the results above together using constructions from one model into the other.

\begin{thm}
\label{thmmain}
Let $(N,*)$ be a rectangular groupoid. Then the graphs $(N,R)$ and $(N,G)$ with $R=\{(a,a*b): a,b \in N\}$ and 
$G=\{(a*b,b): a,b \in N\}$ satisfy property $P_2$.
\end{thm}
Proof: Let $a,b\in N$ be two nodes. Then there is a red edge $(a,a*b) \in R$ and a green edge $(a*b,b)\in G$ so we have at least one red-green path from $a$ to $b$.

Suppose there is a second red-green path from $a$ to $b$, $(a,x) \in R$, $(x,b) \in G$.
 Then there exist some $c,d \in N$ such that $x=a*c$ and $x=d*b$.
 By the rectangularity property, $x=a*b$ so there is no second red-green path and we are done.
 \hfill$\Box$

For any  groupoid we can define such a graph pair, the properties of which will depend upon the properties of the algebra. For instance quasigroups (i.e. the groupoid derived from a Latin square) and only quasigroups will give us two complete graphs.
Commutative idempotent semigroups give us the graphs that are the Hasse diagram of the semilattice order $a\leq b \Leftrightarrow a*b=b$ derived from the operation and the dual order.
A groupoid in general will give us at least one red-green path between any pair of nodes.

\begin{thm}
Let two graphs $(N,R)$ and $(N,G)$ have property $P_2$, 
so for every $a,b\in N$ there is some unique $c \in N$ such that $(a,c)\in R$ and $(c,b)\in G$.
Define $a*b = c$. Then $(N,*)$ is a rectangular groupoid.
\end{thm}
Proof: Suppose $a*b=c*d=x$. Then $(a,a*b)=(a,x) \in R$ and $(x,d)=(c*d,d) \in G$ 
so there is a red-green path from $a$ to $d$ via $x$ and this is unique, so $a*d=x$
\hfill$\Box$

The constructions are exact inverses of one another, so the graph pair derived from the groupoid derived from a graph pair is the same as the original graph pair.

Let's consider a few examples.

\begin{exa}
Take the farmer's market example above with $a=1$. This gives us the array $M$ filled entirely with 1s having property $P_1$, red graph having edges $(x,1)\, \forall x$ and green graph $(1,x)\, \forall x$, the rectangular groupoid with $x*y=1$ for all $x,y$ and the matrices $A$ having all 1s in the first column and zeros elsewhre, $B$ having 1s in the first row and 0s elsewhere such that $AB=J$.

\[
M= \begin{array}{cccc}
1 & 1 & 1 & 1 \\
1 & 1 & 1 & 1 \\
1 & 1 & 1 & 1 \\
1 & 1 & 1 & 1 
 \end{array}, \,\,
A = \left( \begin{array}{cccc}
1 & 0 & 0 & 0 \\
1 & 0 & 0 & 0 \\
1 & 0 & 0 & 0 \\
1 & 0 & 0 & 0 
 \end{array} \right) , \,\,
B = \left( \begin{array}{cccc}
1 & 1 & 1 & 1 \\
0 & 0 & 0 & 0 \\
0 & 0 & 0 & 0 \\
0 & 0 & 0 & 0 
 \end{array} \right) , 
\]
\end{exa}
A somewhat less trivial example 

\begin{exa}
\label{exncg}
Start from the array, graph and matrix as follows:

$
A = \begin{array}{ c c c c }
 1& 1& 2 & 2  \\
 3 & 3 & 4 & 4 \\
 3 & 3& 4 & 4 \\
 1 & 1 & 2 & 2 
\end{array}, 
B = \left( \begin{array}{cccc}
1 & 1 & 0 & 0 \\
0 & 0 & 1 & 1 \\
0 & 0 & 1 & 1 \\
1 & 1 & 0 & 0 \\
 \end{array} \right) , 
(N,E)=$
   \begin{minipage}[]{0.4\textwidth}
      \vspace{0pt}
      \epsfig{file=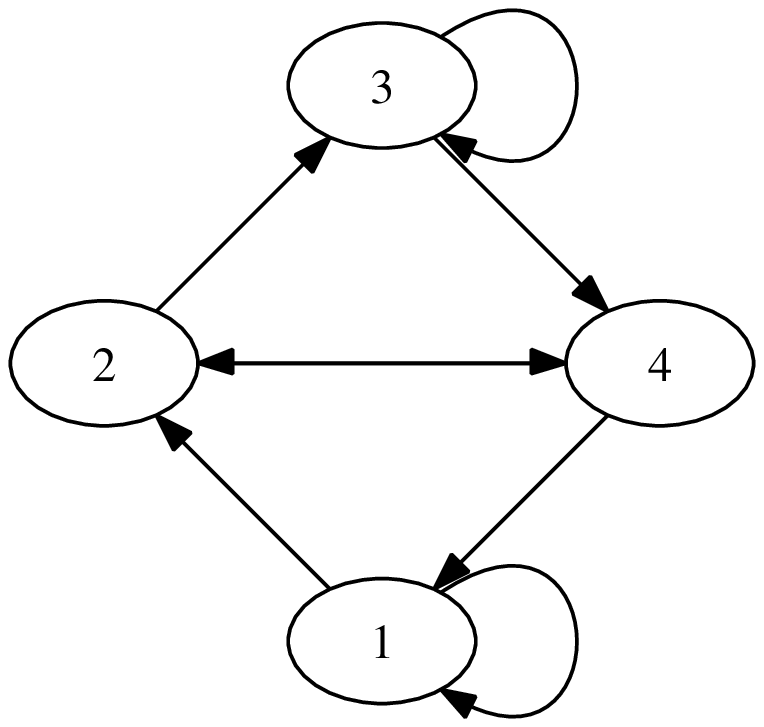, height=30mm}
   \end{minipage}

Then $(N,E),(N,E)$ is a graph pair satisfying $P_2$ corresponding to the array $A$ satisfying $P_1$, the resulting groupoid with Cayley table $A$ and the matrices $BB=J$.
\end{exa}

Note that if $A$ satisfies $P_1$ then so does the transpose $A^T$. This dual structure has a correlate for all the properties above.
\begin{itemize}
\item 
If $R$ is a set of pairs, let $\bar R = \{(b,a): a,b) \in D\}$. Then the dual of a graph pair $(N,R),(N,G)$ is the 
graph pair $(N,\bar G),(N,\bar R)$. A graph pair satisfies $P_2$ iff its dual does.
\item 
The \emph{opposite} groupoid $(N,*)^{opp}$ of a groupoid $(N,*)$ is $(N,+)$ with $a+b=b*a$.  $(N,*)$ is rectangular iff $(N,*)^{opp}$ is.
\item 
Let $A,B$ be $0,1$-matrices. Then $AB=J $ iff $B^T A^T = J$.
\end{itemize}

\section{Special Cases and Related Structures}

A number of special classes of these structures exist and some have been studied previously. In this section we will look at some of these classes, their properties and the way that the various models interrelate.

\subsection{Central groupoids and $UPP_2$ graphs}
\label{seccentral}
In \cite{evans67a} Trevor Evans defined for a set $A$ the groupoid $(A\times A,*)$ with 
\[(a,b)*(c,d) = (b,c)
\]
 These satisfy the equation $(x*y)*(y*z)=y$ and groupoids satisfying this equation are called \emph {central groupoids}. 
 
In \cite{knuth70} Knuth shows that these correspond to $0,1$-matrices $B$ such that $BB=J$ which are equivalent to directed graphs with a unique path of length 2  ($UPP_2$)between all node pairs \cite{lamlint78,mendelsohn70,ryser70}. 
Using the $P_4$ matrix formulation it can be shown that the order of these structures must be a square. The matrices have received special attention, e.g.\ \cite{shader74} showing tight bounds on the possible ranks of the matrices, while circulant matrices have been more specifically investigated \cite{ma84,wu02,wang80,wang82}.
Efforts to exhaustively enumerate small examples(e.g.\ \cite{boykett07,curtisEtAl04} stalled at order $3^2$ with 6 examples until Georg Leander et al, motivated by applications in switching theory, found 3492 examples of order $4^2$ in \cite{leanderEtAl11}.

\subsection{Associativity}
\label{secassoc}

In \cite{mclean54} the class of \emph{rectangular bands} was introduced. A rectangular band $(S,*)$ satisfies the identity $a*b*c=a*c$ as well as associativity and idempotence and are all constructed from two sets $A,B$ with $S=A\times B$ and $(a,b)*(c,d) = (a,d)$.

Let $(A,*)$ be a RG with some $I\subsetneq A$ such that for all $a,b \in A$, $a*b \in I$. Thus the associated $P_1$ array is not full. We call $A$ a  \emph{blow up} of $I$.

As an example, let $(A,*)$ be a rectangular groupoid, $n \not\in A$. Define $n*n=a$, $n*x=a*x$ and $x*n=x*a$ , then $(A\cup\{n\},*)$ is a rectangular groupoid, a blow up of $(A,*)$ by $a$.

The farmers market example above is a blow up of the single element RG $\{a\}$.

\begin{lemma} Let $(A,*)$ be an associative RG. Then $I=\{xy | x,y \in A\}$ is a rectangular band and a subsemigroup of $A$. If $I=A$ then $A$ is a rectangular band, otherwise $A$ is  blow up of $I$.
\end{lemma}
Proof: We write the operation in $A$ as juxtaposition. Let $ab \in I$ then $(ab)(ab) = a(bab)=(aba)b$ so $ (ab)(ab)=ab$ showing that elements of $I$ are idempotent. $I$ consists of all products so it is closed by definition, making it an idempotent subsemigroup. Take $a,b,c \in I$. Then $a(bc)=(ab)c=ac$ so $I$ is  rectangular. If $I\neq A$
then $A$ is a blow up of $I$.
\hfill$\Box$

Owing to the special structure of the rectangular band, there are many blow ups of a rectangular band possible. Let $(S,*)$ be a rectangular band with set sizes $n=\vert A\vert$ and $m=\vert B\vert$. Then a simple counting argument gives $n^{m-1}m^{n-1}(n+m-1)$ extensions not taking into account isomorphism.

Question: Is every blow up of a rectangular band associative? Blow ups constructed as above are associative, but it is not clear that all blow ups of an associative rectangular groupoid are associative.

In a full associative rectangular groupoid, the resulting graph pairs are unions of disjoint isomorphic complete graphs. One graph is $m$ copies of $K_n$ while the other is $n$ copies of $K_m$ with each $K_n$ intersecting each $K_m$ precisely once. 
We can equivalently think of these as two orthogonal partitions of the given node set $N$ of order $nm$. We will see a generalisation of this construction later.

\subsection{Matrix Symmetry}

If we demand a certain higher degree of symmetry in $P_4$, i.e.\  $AB=BA=J$, we obtain another structure.

\begin{thm}
\label{thmscbg}
 $A,B$ are $0,1$-matrices satisfying $P_4$ with the extra symmetrical equation $BA=J$ iff the groupoid $(N,*)$ is a reduct of the algebra $(N,*,+)$ satisfying the equations
 \[(a*b)+(b*c) =b\mbox{ and }(a+b)*(b+c) = b
 \]
 \end{thm}
 Proof: $(\Rightarrow)$ $BA=AB=J$ so we can translate this directly to the graph pair $(N,R),(N,G)$ satisfying $P_2$ (i.e.\ unique red-green path) and the graph pair  $(N,G),(N,R)$ satisfying $P_2$ (i.e. unique green-red path). These give us two rectangular groupoids $(N,*)$ and $(N,+)$.
 If we look at the edges we know that $(a*b,b)$ is a green edge and $(b,b*c)$ is a red edge.
The green-red path from $a*b$ to $b*c$ goes over the node $b$ so $b=(a*b)+(b*c)$ which is the first equation.
 
 The second equation follows from the same argument with the graph pairs reversed.
 
 $(\Leftarrow)$ Suppose we have an algebra $(N,*,+)$ satisfying the two equations. 
 First we show that the groupoids $(N,*)$ and $(N,+)$ are rectangular.
 Let $a,b,c,d \in N$, suppose $a+b=c+d$. By the conditions, we know $b=(a+b)*(b+b)$ and $c=(c+c)*(c+d)$.
 Then
 \beq
 c+b&=&((c+c)*(c+d))+((a+b)*(b+b))\\ &=& ((c+c)*(a+b))+((a+b)*(b+b)) = a+b
 \eeq
 which is the rectangularity property. Similarly we show rectangularity for $(N,*)$.
 
 We can define the four graphs graphs $(N,R_*),(N,G_*),(N,R_+),(N,G_+)$ from these groupoids
 \beq
 R_*&=&\{(a,a*b): a.b \in N\} \\
 R_+&=&\{(a,a+b): a.b \in N\} \\
 G_*&=&\{(a*b,b): a.b \in N\} \\
 G_+&=&\{(a+b,b): a.b \in N\} 
 \eeq
 We will now show that $R_*=G_+$. Let $(a+b,b)\in G_+$, 
 then 
 \beq
 (a+b,(a+b)*(b+c)) = (a+b,b)\in R_*
 \eeq
  so $G_+ \subseteq R_*$. 
 Similarly for all $(a,a*b) \in R_*$, 
 \beq
 ((c*a)+(a*b),(a*b)) = (a,(a*b)) \in G_+
 \eeq
  so $R_* = G_+$. 
 
 Similarly we see that $G_*=R_+$.
 
 Thus we obtain the graph pair $(N,R_*),(N,G_*)$ with the incidence matrices $A,B$ so that $AB=J$. 
 Since the graph pair $(N,R_+),(N,G_+) = (N,G_*),(N,R_*)$ we obtain that $BA=J$ and we are done.
  \hfill$\Box$
 
 This algebraic structure is important in the analysis of reversible one dimensional cellular automata with Welch index not equal to 1\cite{boykett04}.

\subsection{Undirected Graphs}
One can naturally ask when the graph pairs satisfying $P_2$ are undirected, i.e.\ every edge $(a,b)$ has the opposite edge $(b,a)$. The associative case shows that this is possible.

\begin{thm}
Graph pairs satisfying $P_2$ are undirected iff the associated groupoid $(N,*)$ satisfies the equation $(a*b)*(c*a)=a$.
\end{thm}
Proof:
$(\Rightarrow)$ Suppose the graph pair $(N,R),(N,G)$ is undirected and let $(N,*)$ be the associated rectangular groupoid.
We have $(a,a*b) \in R$ so by virtue of the graph being undirected, $(a*b,a)\in R$ too. $(c*a,a) \in G$ implies that $(a,c*a) \in G$ sp we have a red-green path $a*b \rightarrow a \rightarrow a*c$, so $(a*b)*(c*a) = a$ and we are done.

$(\Leftarrow)$  Suppose $(N.*)$ satisfies the equation. Let $(a,a*b)\in R$ be some red edge in the associated graph pair. Then $(a*b,(a*b)*(c*a)) = (a*b,a) \in R$ so $R$ is undirected. Similarly $G$ is undirected and we are done.
\hfill$\Box$

Note that in this case the algebra $(N,*,+)$ with $a+b=b*a$ satisfies the equations given in Theorem \ref{thmscbg}.

This can also be seen directly. The graphs are undirected iff the incidence matrices $A,B$ are symmetrical, i.e.\ $A^T=A$ and $B^T=B$.
Then $BA= B^TA^T = (AB)^T = J^T = J$ so we have the symmetric matrix case from Theorem \ref{thmscbg}.

Note also that it is possible for one graph to be undirected and the other directed, for instance the construction in Section \ref{secpart} below.

\section{Constructions and Reductions}

We investigate several constructions of these structures. First we will look at isotopism as a more general sense of equivalence. Then we will look at substructures, homomorphisms and product constructions.

\subsection{Isotopism as Equivalance}

Given an array $A$ it is clear that reordering the columns or rows or permuting the  entries in the array does not change whether or not the arrays satisfies property $P_1$. The resulting change in the associated groupoid multiplication is called an \emph{isotopy}. An isomorphism is a special type of isotopy.

\begin{defn}
Two groupoids $(A,+)$ and $(B,*)$ are \emph{isotopic} iff $\exists\alpha,\beta,\gamma : B\rightarrow A$ such that
for all $a,b,c \in B$, 
$\alpha(a)+\beta(b) = \gamma(a*b)$.
\end{defn}

The associated graph pair is changed more significantly. Let $(N,R_*),(N,G_*)$ and $(N,R_+),(N,G_+)$ be the graph pairs associated with these two groupoids. The edge $(a,a*b)\in R_*$ is taken to the edge $(\alpha(a),\gamma(a*b))$ by the isotopism, 
$(a*b,b)\in G_*$ is taken to $(\gamma(a*b),\beta(b)) \in G_+$. 
That is, $(a,b)\in R_*$ is mapped to $(\alpha(a),\gamma(b)) \in R_+$, a more significant change.

We see here that isotopies indicate several distinct degrees of ``sameness.'' 
An isotopy of arrays satisfying $P_1$ gives us something essentially the same, the same applied to the rectangular groupoid is less identical. Applying an isotopy, the associated graph pair is definitely different: for instance loop edges may arise or disappear.

A \emph{transversal} of an array satisfying $P_1$ of order n is a set of n cells with the property that one cell lies in each row, one in each column, and one contains each symbol. 

\begin{thm}
An array satisfying $P_1$ has a transversal iff the associated groupoid has an idempotent isotope.
\end{thm}
Proof: Let $n$ be the size of the array $M$. 
Let  $A=\{1,\ldots,n\}$ and  $(A,*)$ be the associated rectangular groupoid, i.e.\ $a*b$ is the entry in row $a$ and column $b$. 

$(\Rightarrow)$: Let the vectors $v,w \in \{1,\ldots,n\}^{ \{1,\ldots,n\}}$ have $v(i)$ being the row where $i$ occurs in the transversal, $w(i)$ be the column where $i$ appears in the transversal. 
Then the mappings $\alpha:i\mapsto v(i)$ and $\beta: i \mapsto w(i)$ are permutations of $A$. 
The isotopy $(\alpha^{-1},\beta^{-1},id)$ maps 
$(A,*)$ to $(A,+)$ with $a+b = \alpha(a)*\beta(b)$. 
Now $\alpha(i)$ is the row where $i$ appears in the transversal, $\beta(i)$ is the column where $i$ appears. 
So entry $(\alpha(a),\beta(a))$ is $a$, thus $a+a = \alpha(a)*\beta(a) = a$.

$(\Leftarrow)$:  Suppose the rectangular groupoid $(A,+)$ is idempotent and isotopic to $(A,*)$ associated with the array $M$ by the isotopy $(\alpha,\beta,\gamma)$. That is, $\gamma(a*b) = \alpha(a)+\beta(b)$.
Let $T = \{(\alpha^{-1}(i),\beta^{-1}(i)); i \in  \{1,\ldots,n\}\}$.
The $(i,j)$ entry in $M$ is $i*j=  \gamma^{-1}( \alpha(i)+\beta(j))$, so the $(\alpha^{-1}(i),\beta^{-1}(i))$ entry in $M$ is
\[
\alpha^{-1}(i) * \beta^{-1}(i) = \gamma^{-1}(\alpha(\alpha^{-1}(i))+\beta(\beta^{-1}(i))) = \gamma(i+i) = \gamma(i)
\]
Because $\gamma$ is a permutation, this means that $T$ is a transversal of $M$ and we are done.
\hfill$\Box$

In the case of matrix symmetric rectangular groupoids, we know that every example is isotopic to a unique idempotent matrix symmetric rectangular groupoid \cite{boykett04}.
The following question then arises: can a rectangular groupoid be isotopic to two nonisomorphic idempotent rectangular groupoids?
The answer here is no. The following examples have been found from an exhaustive listing generated by Mace \cite{prover9-mace4}.

\[
\begin{array}{c|ccccc}
* & 0 & 1 & 2  & 3 & 4\\
\hline
0 & 0 & 0 & 0 & 0 & 0 \\
1 & 1 & 1 & 2 & 2 & 1 \\
2 & 1 & 1 & 2 & 2 & 1 \\
3 & 3 & 4 & 3 & 3 & 4 \\
4 & 3 & 4 & 3 & 3 & 4
\end{array}
\,\,\,\,\,\,\,\,\,\,\,\,
\begin{array}{c|ccccc}
+ & 0 & 1 & 2  & 3 & 4\\
\hline
0 & 0 & 0 & 0 & 0 & 0 \\
1 & 1 & 1 & 2 & 2 & 2 \\
2 & 1 & 1 & 2 & 2 & 2 \\
3 & 3 & 3 & 4 & 3 & 4 \\
4 &3 & 3 & 4 & 3 & 4 
\end{array}
\]

The two rectangular groupoids are not isomorphic but are isotopic by the column permutation $\beta=(0 3 1 2)$ and entry permutation $\gamma = (1 2)$. Thus the idempotent examples cannot be used as representatives of each isotopy class as in the matrix symmetric case.

Question: is there a subvariety $V$ of the quasivariety of rectangular groupoids such that every rectangular groupoid is isotopic to exactly one in $V$? It is not sufficient to restrict ourselves to full rectangular groupoids, as we see by these examples (both are full). The class of matrix symmetric rectangular groupoids is too small, as every isotope of a matrix symmetric rectangular groupoid is matrix symmetric.

\subsection{Substructures, Products and Homomorphic images}

There are a number of methods available to take structures and combine them to obtain new ones. Some of the classical methods are to take direct products, homomorphic images and substructures. These are most easily applied to the algebraic formulation as groupoids. 



Because the class of rectangular groupoids has been written with a defining quasiidentity (\ref{exrect}), we know that the class forms a quasivariety and thus is closed under the taking of subalgebras and direct products. 

However we can demonstrate that the quasivariety of rectangular groupoids is particularly badly behaved.

\begin{thm}
The smallest variety containing the rectangular groupoids is the  variety of all groupoids.
\end{thm}
Proof: We demonstrate this by showing that for all groupoids $(G,*)$ there is a rectangular groupoid with $(G,*)$ as a homomorphic image.
Let $(G,*)$ be a groupoid. Define an operation $+$ on $G\times G$
by $(a,b)+(c,d) = (a*c,c)$. First we show that $(G\times G,+)$ is rectangular, then we will show it has $G$ as a homomorphic image.

Suppose $(a,b)+(c,d) = (\bar a,\bar b)+(\bar c, \bar d)$. Then $c=\bar c$ and $a*c= \bar a * \bar c= \bar a * c$.
Thus $(a,b)+(\bar c, \bar d) = (a*\bar c,\bar c) =(a*c,c)= (a,b)+(c,d)$ and similarly 
$(\bar a,\bar b)+(c,d) = (a,b)+(c,d)$ so we see rectangularity of $(G\times G,+)$.

The map $\alpha:G\times G \rightarrow G$, $(a,b)\mapsto a$ is an epimorphism so $(G,*)$ is a homomorphic image of the rectangular groupoid $(G\times G,+)$ so variety generated by rectangular groupoids is all groupoids.
\hfill$\Box$

Thus there are no nontrivial equations satisfied by all rectangular groupoids.

One can see this less clearly but more easily using the associated graph pair. In the homomorphic image of such a graph pair, we will still have the condition that at least one red-green path exists between each pair of nodes, but we will not be able to claim that this path is unique, as the graph homomorphism may map the end points of two paths together but not the middle nodes.

Many subclasses of rectangular groupoids are varieties as we have seen above. One of the most natural subclasses  are the idempotent rectangular groupoids.  The following example shows that these are also not closed under taking homomorphic images. We take the congruence with partition $1,2,3,4\vert n$ to form the homomorphism.

\begin{center}
\[
\begin{array}{ccccc}
1 & 1 & 3 & 3 & 3 \\
2 & 2 & 4 & 4 & 2 \\
1 & 1 & 3 & 3 & 3 \\
2 & 2 & 4 & 4 & 2 \\
1 & 1 & 4 & 4 & n 
 \end{array}
 \xrightarrow{{1234\vert n}}
 \begin{array}{cc}
1 & 1  \\
1 &  n 
\end{array}
\]
\end{center}

However the situation is not  as with general rectangular groupoids.

\begin{thm}
The variety generated by idempotent rectangular groupoids is a proper subvariety of the idempotent groupoids.
\end{thm}
Proof: Let $(N,*)$ be an idempotent rectangular groupoid, $a,b \in N$. Then $(a*b)*(a*b)=a*b$ by idempotence, so
$a*(a*b) = (a*b)*b=a*b$. Since these equations hold for all  idempotent rectangular groupoids they also hold for the generated variety $\V$. The groupoid $(\{0,1,2\},*)$ defined by the table
\[
\begin{array}{c|ccc}
* & 0 & 1 & 2 \\
\hline
0 & 0 & 2 & 2 \\
1 & 0 & 1 & 2 \\
2 & 0 & 1 & 2
\end{array}
\]
is idempotent but does not satisfy the equation because $(0*1)*1=2*1=1$ but $0*1=2$. Thus this groupoid is not in $\V$ so the $\V$ is properly contained in the variety of all idempotent groupoids.
\hfill$\Box$

\subsection{Partition Construction Technique}
\label{secpart}

Let $\Pi$ be a partition of $N$ and for every part $\pi \in \Pi$ let $\theta_\pi$ be a partition of $N$ with $\pi$ a transversal of $\theta_\pi$. Let $(N,R)$ be the graph formed by union of complete graphs on each part $\pi$. 
Let $G=\{(a,b) \in \theta_\pi \mbox{ with } a \in \pi\}$. Then $(N,R),(N,G)$ is a graph pair satisfying $P_2$. We call such a structure \emph{partitioned}.
This generalises a construction suggested by Tim Penttila for matrix symmetric rectangular groupoids.

\begin{thm}
Let $(N,R),(N,G)$ be a graph pair satisfying $P_2$. Then the following are equivalent:
\begin{enumerate}
\item $(N,R)$ is a union of cliques and $(N,G)$ has loops on each node
\item  $(N,R),(N,G)$ is partitioned
\item the associated groupoid satisfies the equations $a*a=a$ and $(a*b)*c=a*c$
\end{enumerate}
\end{thm}
Proof:
$(2)\Rightarrow(1)$ follows from the construction.

$(1)\Rightarrow (3)$: 
Let $a,b,c \in N$. There is a red loop edge on $a$ and a green loop edge on $a$ so the path from $a$ to $a$ goes through $a$ so $a*a=a$.
The nodes $a$ and $a*b$ are in the same red clique, as are the nodes $a*b$ and $(a*b)*c$, so all three are in the same clique so there is a red edge from $a$ to $(a*b)*c$. Because there is a green edge from $(a*b)*c$ to $c$ then there is a red-green path from $a$ to $c$ via $(a*b)*c$ so $a*c=(a*b)*c$.

$(3)\Rightarrow (1)$: 
Because $a*a=a$ we have a red and a green loop edge on each node.
Suppose $(a,b),(b,c) \in R$, that is there exist $n,m \in N$ such that $b=a*n$ and $c=b*m=(a*n)*m$. But then $a*m=(a*n)*m$ by condition $(3)$, so  $(a,c)\in R$. Thus $R$ is reflexive and transitive. Now $b*a = (a*n)*a = a*a = a$ so $(a,b)\in R \Rightarrow (b,a)\in R$ so $R$ is symmetric and thus an equivalence relation, so $(N,R)$ is a union of cliques.

$(1)\Rightarrow (2)$: 
Let $\Pi$ be the partition induced by the cliques in $R$. Let $\pi \in \Pi$ be one part. Suppose there exists $a,b \in \pi$ and $c\in N$ with $(a,c),(b,c) \in G$. Then because $(a,a),(a,b)\in R$ there exists two red-green paths from $a$ to $c$ which is a contradiction. So the green edges leaving $\pi$ partition $N$. Call this partition $\theta_\pi$. Then we are done.
\hfill$\Box$

Note that in a $P_2$ graph pair $(N,R),(N,G)$  is a union of cliques iff the above theorem applies in the dual graph pair. 
If the dual of a graph pair is partitioned we say that the graph pair is dually partitioned. 
The following result is immediate.

\begin{cor}
Let $(N,R),(N,G)$ be a graph pair satisfying $P_2$. Then the following are equivalent:
\begin{enumerate}
\item $(N,G)$ is a union of cliques and $(N,R)$ has loops on each node
\item  $(N,R),(N,G)$ is dually partitioned
\item the associated groupoid satisfies the equations $a*a=a$ and $a*(b*c)=a*c$
\end{enumerate}
\end{cor}

We have seen the following result above in a different form, the two partitions are generated by the sets $A,B$ that give the rectangular band $A\times B$.

\begin{lemma}
A graph pair satisfying $P_2$ is partitioned and dually partitioned iff it is associative.
\end{lemma}

We can create such examples from groups.
Let $\Gamma$ be a group, $H\leq \Gamma$ a subgroup and $1\in T\subset \Gamma$ a set of left coset representatives of $H$ in $\Gamma$. Then the left cosets of $H$ form  a partition and for each part $aH$ the partition from the equivalence relation  $\{(ah,ath) :h\in H,\, t\in T\}$ has $aH$ as a transversal.

In this case the red graph is a collection of cliques and the green 
graph is the Cayley graph with node set $\Gamma$ generated by $T$.

This idea can be extended to any set factorisation of a group $\Gamma$ into two subsets $H,K \subset \Gamma$ with $HK=\Gamma$ and $\vert H \vert \vert K \vert = \vert \Gamma \vert$. Then every element of $\Gamma$ has a unique representation as $hk$ for some $h\in H, k \in K$ and the Cayley graphs on $\Gamma$ generated by $H$ and $K$ form a graph pair with $P_2$.

The following result follows in a similar way to the recognition of difference families in BIBDs \cite{beth85}.

\begin{thm}
A $P_2$ graph pair has a regular automorphism group iff it is two Cayley graphs as described above.
\end{thm}
Proof: $(\Rightarrow)$ 
Let $\Gamma$ be the regular automorphism group acting on the left. Identify $N$ and $\Gamma$ so $\Gamma$ acts on itself by left multiplication. Let $H\subseteq \Gamma$ be the set of red neighbours of the identity $1\in \Gamma$, $K\subseteq \Gamma$ be the set of green neighbours of the identity.

We claim that $E_H= \{(a,ah): a \in \Gamma, h\in H\}$ is the set of red edges. Let $(a,b)\in R$ be a red edge. Then we apply the automorphism $a^{-1}$ to see the edge $(1,a^{-1}b)$ so $a^{-1}b\in H$ so $(a,b) \in E_H$, $R\subseteq E_H$.
Likewise all members of $E_H$ are images of a red edge starting from the identity so $E_H \subseteq R$ and we are done. Similarly all green edges are generated by $K$.

$(\Leftarrow)$ : The group $\Gamma$ acting by left multiplication takes edges to edges, $\alpha (a,ak) = ((\alpha a,(\alpha a) k)$ and is a regular automorphism group of both graphs.
\hfill$\Box$

\subsection{Combining Rectangular Groupoids}

Given two rectangular groupoids, there are a number of ways of combining them to create a new rectangular groupoid.

Let $A,B$ be two rectangular groupoids, $A \cap B = \emptyset$ and $f:A\rightarrow B$, $g:B \rightarrow A$ two mappings. We define a new rectangular groupoid on $A \cup B$ with
\[
\begin{array}{|c|c|}
\hline
 & \mbox{columns}\\
A & \mbox{of A}\\
\hline
\mbox{columns } & B  \\
\mbox{of B} &   \\
\hline
\end{array}
\hspace{12mm}
x*y = \left\{ 
\begin{array}{cl}
x *_A y &\mbox{ if } x,y \in A \\
x *_A g(y)& \mbox{ if } x \in A,\, y \in B \\
f(x) *_B y &\mbox{ if } x \in B,\, y \in A \\
x *_B y &\mbox{ if } x,y \in B
\end{array}
\right .
\]
If we look at the array that arises, we make a block diagonal new array with A and B on the diagonal. We copy columns from the A section into the top right block, columns from $B$ into the bottom left block. We thus add no new pairs of elements appearing in the same row together, the columns receive new pairs from $A\times B$ which do not appear in any rows. So the new array satisfies the conditions of $P_1$ if the  starting arrays $A$ and $B$ do.

We call this a \emph{left split extension} because of the way the left side of products in $A \cup B$ define where the product lies. Similarly we can define a \emph{right split extension} by placing rows of $A$ in the bottom left block and rows of $B$ in the top right block.

We saw an example of this in Example \ref{exncg} where the array is a right split extension of the two associative rectangular groupoids $\{1,2\}$ and $\{3,4\}$. 

Another extension is made as follows. Given an rectangular groupoid $A$ and an element $a \in A$ we create a new element $n \not\in A$, define a new array on $A\cup \{n\}$ with $x*n=x*a$, $n*n=n$ and 
\[
 n*x  = \left\{ 
 \begin{array}{cl}
 a*x &\mbox{ if }a*x \neq a \\
 n& \mbox{ otherwise}
 \end{array}
 \right.
 \]
 We call this the \emph{left extension} of $A$ by $a$. Similarly we define the \emph{right extension} of $A$ by $a$.
 
 Investigating an exhaustive list of all small examples, we see that almost all examples are obtained from a smaller one by one of these extensions. The smallest nonexample is the 5 element example  shown above as a counterexample to the homomorphic closure of idempotent rectangular groupoids.
 
Question: If a class of rectangular groupoids are closed under homomorphisms i.e.\ all homomorphic images of them are in the class, then the split extensions of them and the one element extensions of them are also in the class.
Alternatively, if $A,B$ are rectangular groupoids such that all homomorphic images of them are also rectangular groupoids, then all split and left/right extensions of them also have the property that all homomorphic images are rectangular groupoids.

\section{The common root of Rectangularity and Latinicity}

We introduced these arrays as some kind of opposite of Latin squares. Both concepts can be generalised in the sense that we can talk about incomplete arrays that do not break the requirements of the given structures.

Let $M$ be an $n\times n$ array partially filled with entries from $\{1,\ldots,n\}$.
We say that $M$ is a \emph{partial Latin square} if each row and column contains at most one copy of each element.

$M$ is a  \emph{partial $P_1$-array} if it satisfies $P_1$.
By analogy to Theorem \ref{thmmain} we can say that a partial $P_1$ array corresponds to a graph pair with at most one red-green path between any set of nodes.

A partial Latin square has \emph{Blackburn property} \cite{wanless04}, derived from the construction of perfect hash families \cite{blackburn00}, if whenever the cells $(i,j)$ and $(k,l)$ are occupied by the same symbol, the opposite corners $(i,l)$ and $(k,j)$ are empty. 

\begin{thm}
A partial Latin square  that is also a partial $P_1$-array has the Blackburn property.
\end{thm}
Proof:
Let $M$ be such an array.
Suppose the cells  $(i,j)$ and $(k,l)$ are occupied by the same symbol $a$ and the cell $(i,l)$ is occupied with the symbol $b$. Then the pair $(a,b)$ appears in row $i$ and column $l$  which contradicts $P_1$ unless $a=b$. But then we have two occurences of $a$ in row $i$ and column $l$ which contradicts the Latin square property. So the cell $(i,l)$ is empty, as is $(k,j)$ and we have shown the Blackburn property.
\hfill$\Box$

Unfortunately not all partial Latin squares with the Blackburn property satisfy $P_1$, as demonstrated by
\[
\begin{array}{ccc}
\cdot  & a &  c  \\
a & \cdot   & b \\
c & b & \cdot 
 \end{array}
\]

\section{Conclusions}

We introduced several combinatoric structures and showed that these are all closely related. Several special cases have been investigated previously. We developed connections between these.
While the ideas here are somehow maximally different to those of Latin squares, there is a common core around the idea of the Blackburn property.

The idea of rectangular groupoids can be extended to $n$-ary functions. We say a function $f:A^n\rightarrow A$ is rectangular when  $f(a_1,\ldots,a_n)=f(b_1,\ldots,b_n) \Rightarrow  \forall i\, f(a_1,\ldots,b_i,\ldots,a_n)=f(b_1,\ldots,b_n)$. 
It has been found \cite{boykettkaritaati} that such functions allow a certain amount of ``physical'' behaviour (conservation laws) in one dimensional cellular automata. Related ideas are also known in circuit theory \cite{poeschel92}, their algebras being a special case of $n$-ary rectangularity.

One of the main problems here is that  there are far too many examples. Thus our attention is focussed upon developing descriptions that allow us to investigate a smaller but still important collection of examples, for instance idempotent rectangular groupoids or the various varieties that were introduced above.

\section{Acknowledgements}

The core of this work was developed  on a walking weekend in the comfortable ``Stube'' of  a mountain hut on the Wurzeralm in the Totes Gebirge, and I am most grateful to my hosts. 
Details have been worked out within project P19463 of the Austrian FWF  as well as the Fractured project sponsored by the city of Linz.

\bibliographystyle{elsarticle-num}
\bibliography{rectgpoid}



\end{document}